\documentclass[11pt]{article}
 \usepackage{graphicx}
 \usepackage{amssymb, amsmath,amsthm}
 \usepackage{epstopdf}
\newtheorem{thm}{{\bf Theorem}}[section]
 \newtheorem{lem}[thm]{{\bf Lemma}}
 
 \newtheorem{prop}[thm]{{\bf Proposition}}
 
 \newtheorem{claim}[thm]{Claim} 
 
 \newtheorem{ex}[thm]{Example}

\numberwithin{equation}{section}

\begin{document}
\title{An asymptotic behavior of the dilatation for a family of pseudo-Anosov braids}
\author{Eiko Kin\footnote{The first author is partially supported by Grant-in-Aid for Young Scientists (B) (No. 17740094), 
The Ministry of Education, Culture, Sports, Science and Technology, Japan} \hspace{1mm} and Mitsuhiko Takasawa}
\maketitle

Abstract. 
The dilatation of a pseudo-Anosov braid  is a conjugacy invariant. 
In this  paper, we study the dilatation of a special  family of pseudo-Anosov braids. 
We prove  an inductive formula to compute their dilatation, a monotonicity and an asymptotic behavior of the dilatation for this family of braids. 
We also  give an example of a family of pseudo-Anosov braids with arbitrarily small dilatation such that 
the mapping torus obtained from such braid  has $2$ cusps and has an arbitrarily large   volume. 
\medskip

\medskip


\section{Introduction}

 Let $\Sigma = \Sigma_{g,p}$ be an orientable surface of genus $g$ with $p$ punctures, and let $\mathcal{M}(\Sigma)$ be the mapping class group of $\Sigma$. 
 The elements of $\mathcal{M}(\Sigma)$, called {\it mapping classes}, are classified into $3$ types: periodic, reducible and pseudo-Anosov \cite{Thu}. 
 For a pseudo-Anosov mapping class $\phi$, the dilatation $\lambda(\phi)$ is an algebraic integer strictly greater than $1$.  
The dilatation of a pseudo-Anosov mapping class  is a conjugacy invariant.

 Let $D_n$ be an  $n$-punctured closed disk.  
 The mapping class group $\mathcal{M}(D_n)$ of $D_n$  is isomorphic to a subgroup of $\mathcal{M}(\Sigma_{0,n+1})$. 
 There is a natural surjective homomorphism 
 $$\Gamma: B_n \rightarrow \mathcal{M}(D_n)$$ 
 from the $n$-braid group $B_n$ to  the mapping class group $\mathcal{M}(D_n)$ \cite{Bir}. 
We say that a braid $\beta \in B_n$ is {\it pseudo-Anosov} if $\Gamma(\beta) $ is pseudo-Anosov, and if this is the case 
the dilatation $\lambda(\beta)$ of $\beta$  is defined equal to $\lambda(\Gamma(\beta))$. 
Henceforth, we shall abbreviate $\lq$pseudo-Anosov' to $\lq$pA'.

\begin{figure}[htbp]
\begin{center}
\includegraphics[width=4.5in]{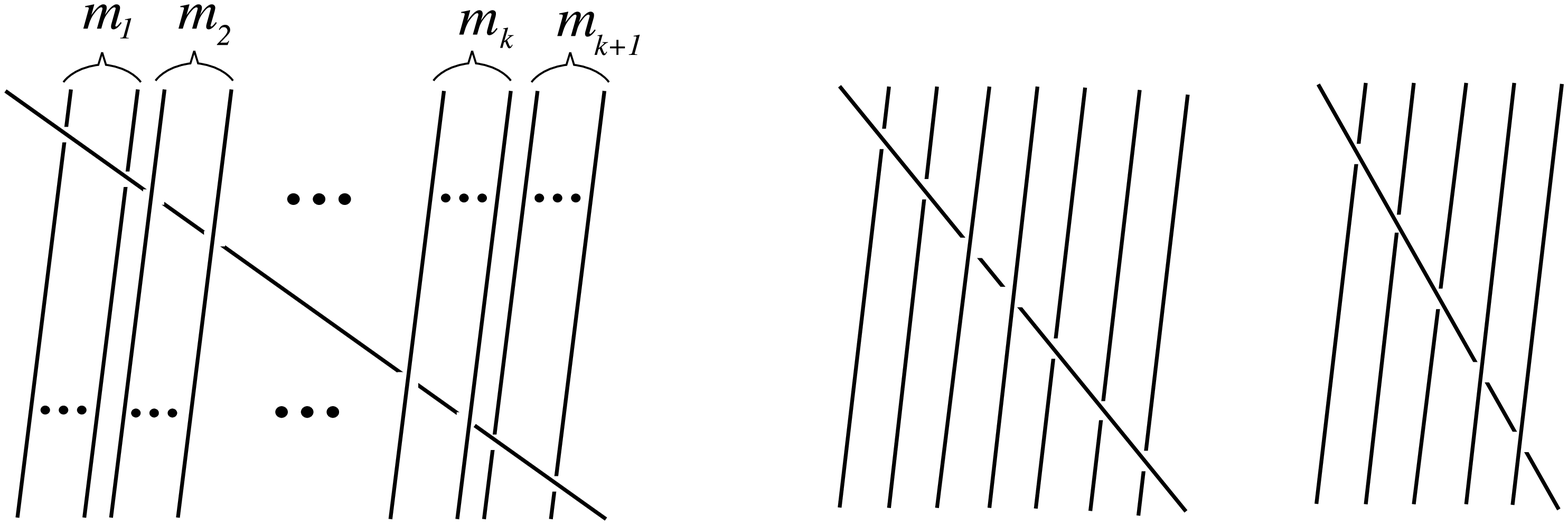}
\caption{(left) $\beta_{(m_1,m_2, \cdots, m_{k+1})}$, (center) $\beta_{(2,2,3)}$, (right) $\beta_{(3,2)}$.}
\label{fig_general_braid}
\end{center}
\end{figure}

We now introduce a family of braids. 
Let $\beta_{(m_1,m_2, \cdots, m_{k+1})}$   be the braid as depicted  in Figure~\ref{fig_general_braid}, 
for each integer $k \ge 1$ and each integer $m_i \ge 1$. 
These are all pA (Proposition~\ref{prop_pseudo-Anosov}). 
We will prove a  monotonicity, an inductive formula to compute their dilatation  and an asymptotic behavior  of the dilatation for this family of braids.

\begin{prop} 
(Monotonicity)
\label{prop_decreasing-example}
For each integer $i$ with  $1 \le i \le k+1$, we have 
$$\lambda(\beta_{(m_1, \cdots, m_i, \cdots, m_{k+1})}) > 
\lambda(\beta_{(m_1, \cdots, m_i+1, \cdots, m_{k+1})}).$$
Hence if $m_i \le m_i'$ for  each $i$, then 
$\lambda(\beta_{(m_1, \cdots, m_{k+1})}) \ge \lambda(\beta_{(m'_1, \cdots, m'_{k+1})})$.  
\end{prop}

For an integral polynomial $f(t)$ of degree $d$, the {\it reciprocal} of $f(t)$, denoted by $f_*(t)$, is $t^d f(1/t)$.

\begin{thm} 
(Inductive formula)
 \label{thm_recursive-example}
 The dilatation of the pA braid $\beta_{(m_1,  \cdots, m_{k+1})}$ is the largest root of the polynomial 
$$ t^{m_{k+1}} R_{(m_1, \cdots, m_k)}(t)+ (-1)^{k+1} {R_{(m_1, \cdots, m_k)}}_*(t),$$
where $R_{(m_1, \cdots, m_i)}(t)$ is given inductively as follows: 
\begin{eqnarray*}
R_{(m_1)}(t) &=& t^{m_1+1} (t-1) -2t, \ \mbox{and}
\\
R_{(m_1, \cdots, m_i)}(t)&=& t^{m_i} (t-1) R_{(m_1, \cdots, m_{i-1})}(t)+ (-1)^i 2t {R_{(m_1, \cdots, m_{i-1})}}_*(t)\ \mbox{for\ } 2 \le i \le k. 
\end{eqnarray*}
 \end{thm}

\begin{thm} 
(Asymptotic behavior)
 \label{thm_lim-example} 
 We have 
\medskip
\\
 (1) 
$\displaystyle\lim_{m_1, \cdots, m_{k+1} \to \infty} \lambda(\beta_{(m_1, \cdots, m_{k+1})}) = 1$ and 
 \medskip
 \\
(2)\
$ \displaystyle\lim_{m_i, m_{i+1}, \cdots, m_{k+1} \to \infty} \lambda(\beta_{(m_1,\cdots, m_{k+1})}) = \lambda(R_{(m_1, \cdots, m_{i-1})}(t)) >1$ 
for $i \ge 2$, where $\lambda(f(t))$ denotes the maximal absolute value of the roots of $f(t)$. 
\end{thm}

For a pA braid $\beta$, let $\phi$ be the pA mapping class $\Gamma(\beta)$.  
The dilatation $\lambda(\phi)$ can be  computed as follows. 
A smooth graph $\tau$, called a {\it train track} and a smooth graph map  $\widehat{\phi}: \tau \rightarrow \tau$ are associated with $\phi$. 
The edges of $\tau$ are classified into {\it real} edges and {\it infinitesimal} edges, and 
the {\it transition matrix} $M_{\mathrm{real}}(\widehat{\phi})$ with respect to real edges can be defined. 
Then the dilatation $\lambda(\phi)$ equals the spectral radius of $M_{\mathrm{real}}(\widehat{\phi})$. 
For more details, see Section~\ref{subsection_traintrack}.

\begin{figure}[htbp]
\begin{center}
\includegraphics[width=6in]{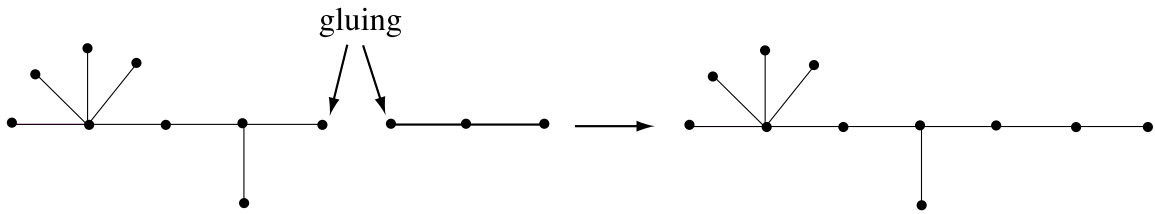}
\caption{$\mathcal{Q}_{(4,2,1)}$ (right) is obtained by gluing $\mathcal{Q}_{(4,2)}$ (left) and another tree (center).}
\label{fig_ex-combined}
\end{center}
\end{figure}

\begin{figure}[htbp]
\begin{center}
\includegraphics[width=6in]{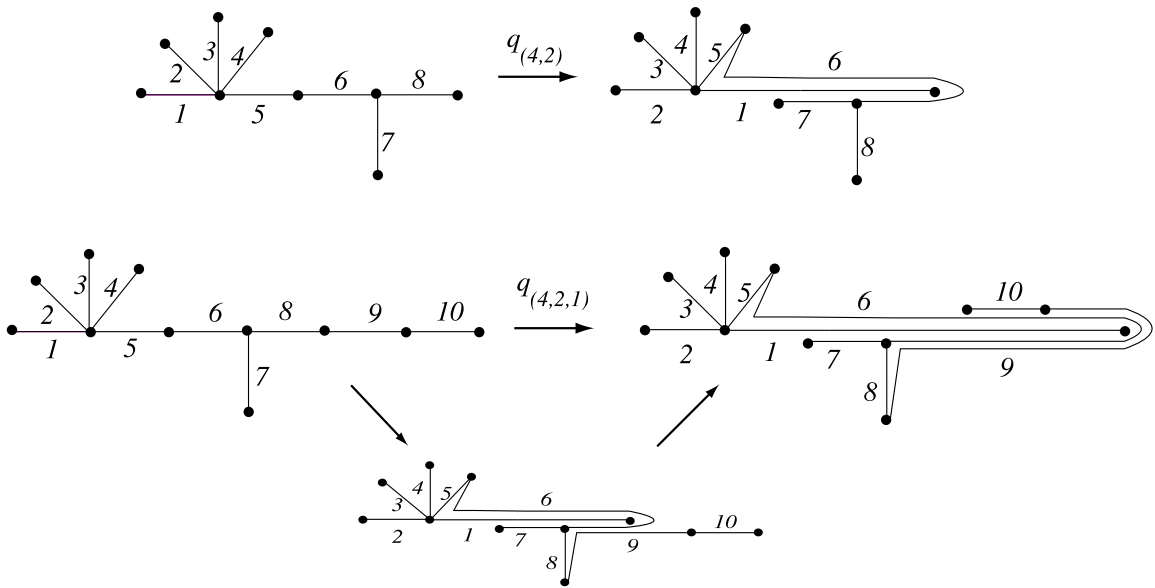}
\caption{(top) $q_{(4,2)}$, (bottom) $q_{(4,2,1)}$.}
\label{fig_ex-combined-map}
\end{center}
\end{figure}

For the computation of the dilatation of the braid $\beta_{(m_1, \cdots, m_{k+1})}$, we introduce {\it combined trees} and {\it combined tree maps} 
in Section~\ref{section_combined}. 
For a given $(m_1, \cdots, m_{k+1})$, one can obtain the combined tree $\mathcal{Q}_{(m_1, \cdots, m_{k+1})}$ and the combined tree map 
$q_{(m_1, \cdots, m_{k+1})}$ inductively. 
For example, for $(m_1,m_2,m_3)= (4,2,1)$,  the combined tree $\mathcal{Q}_{(m_1,m_2,m_3)}$, depicted in Figure~\ref{fig_ex-combined},  is obtained 
by gluing the combined tree $\mathcal{Q}_{(m_1,m_2)}$ and another tree which depends $m_3$. 
The combined tree map $q_{(m_1,m_2,m_3)}$, as shown in Figure~\ref{fig_ex-combined-map},  is defined by the composition of an extension of the combined tree map $q_{(m_1,m_2)}$ and another tree map which depends on $m_3$.

By the proof of Proposition~\ref{prop_pseudo-Anosov}, it turns out that 
the spectral  radius of the  transition matrix $M(q_{(m_1, \cdots, m_{k+1})})$ obtained from $q_{(m_1, \cdots, m_{k+1})}$ equals 
that of $M_{\mathrm{real}}(\widehat{\phi})$, where $\phi= \Gamma(\beta_{(m_1, \cdots, m_{k+1})})$, 
that is the spectral radius of $M(q_{(m_1, \cdots, m_{k+1})})$ equals the dilatation  $\lambda(\beta_{(m_1, \cdots, m_{k+1})})$. 
Proposition~\ref{prop_decreasing-example} and Theorems~\ref{thm_recursive-example}, \ref{thm_lim-example} 
will be shown by using the properties of combined tree maps.

In the final part, we will consider the two invariants of pA mapping classes, the dilatation and the volume. 
 Choosing any representative $f: \Sigma \rightarrow \Sigma$ of  a mapping class $\phi$, we form the mapping torus 
 $${\Bbb T}(\phi)= \Sigma \times [0,1]/ \sim, $$ 
 where $\sim$ identifies $(x,0)$ with $(f(x),1)$. 
 A mapping class $\phi$ is pA if and only if ${\Bbb T}(\phi)$ admits a complete hyperbolic structure of finite volume \cite{Ota}. 
 Since such a structure is unique up to isometry, 
 it makes sense to speak of  the {\it volume} $\mathrm{vol}(\phi)$ of $\phi$,  the hyperbolic volume of ${\Bbb T}(\phi)$. 
 For a pA braid $\beta$, we define the volume $\mathrm{vol}(\beta)$ as equal to $\mathrm{vol}(\Gamma(\beta))$, 
 the volume of the mapping torus ${\Bbb T}(\Gamma(\beta))$. 
  
  Theorem~\ref{thm_lim-example}(1) tells us that dilatation of braids can be  arbitrarily small. 
 We consider what  happen for the volume of a family of pseudo-Anosov mapping classes  whose dilatation is arbitrarily small. 
It is not hard to see the following. 
 
\begin{prop}
\label{prop_trivial-construction}
There exists a  family of pA mapping classes $\phi_n$ of $\mathcal{M}(D_n)$ such that 
$$\lim_{n \to \infty} \lambda(\phi_n) = 1\ \mbox{and}\  \lim_{n \to \infty} \mathrm{vol} (\phi_n) = \infty$$ 
and such that the number of the cusps of the mapping torus ${\Bbb T} (\phi_n)$ goes to $\infty$ as $n$ goes to $\infty$. 
\end{prop}
 
\noindent
Proposition~\ref{prop_trivial-construction} is not so surprising, because the volume of each cusp is bounded below uniformly. 
We show the following.

\begin{prop}
\label{prop_volume-dilatation}
There exist a  family of pA mapping classes $\phi_n$ of $\mathcal{M}(D_n)$ such that 
$$\lim_{n \to \infty} \lambda(\phi_n) = 1\ \mbox{and}\  \lim_{n \to \infty} \mathrm{vol} (\phi_n) = \infty$$ 
and such that the number of the cusps of the mapping torus ${\Bbb T}(\phi_n)$ is $2$ for each $n$. 
\end{prop}
 
\noindent
Proposition~\ref{prop_volume-dilatation} is a corollary of the following theorem. 
 
 \begin{thm}
\label{thm_dil-vol}
For any real number $\lambda>1$ and any real number $v >0$, 
 there exist  an integer $k\ge 1$ and an integer $m \ge 1$ such that  
 for any integer $m_i \ge m$ with $1 \le i \le k+1$, we have 
 \begin{center}
 $\lambda(\beta_{(m_1, \cdots, m_{k+1})}) < \lambda$ \hspace{2mm}and\hspace{2mm} 
 $\mathrm{vol}(\beta_{(m_1,  \cdots, m_{k+1})}) >v$. 
 \end{center}
\end{thm}

\begin{figure}[htbp]
\begin{center}
\includegraphics[width=2in]{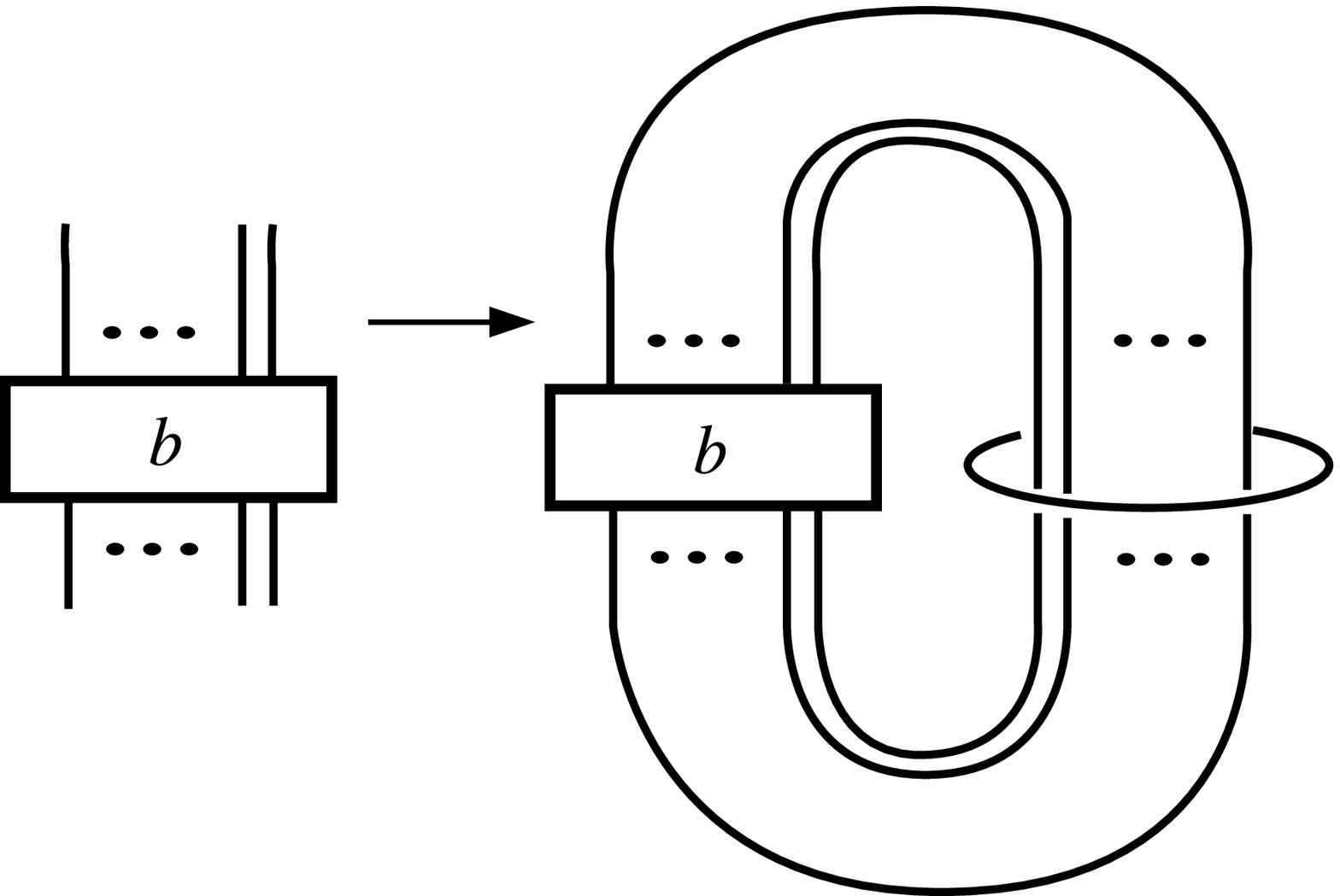}
\caption{link $\overline{b}$.} 
\label{fig_axis}
\end{center}
\end{figure}

\noindent
 Here we note that for a braid $b$, the mapping torus ${\Bbb T}(\Gamma(b))$ is homeomorphic to the link complement $S^3 \setminus \overline{b}$ 
 in the $3$ sphere $S^3$, where $\overline{b}$ is a union of the closed braid of $b$ and the braid axis (Figure~\ref{fig_axis}). 
 When  $b$ is a braid $\beta_{(m_1, \cdots, m_{k+1})}$, the link $\overline{b}$ has $2$ components, and hence 
 the number of cusps of ${\Bbb T}(\Gamma(b))$ is $2$.

\section{Preliminaries}

A homeomorphism $\Phi: \Sigma \rightarrow \Sigma$ is {\it pseudo-Anosov} ({\it pA})  if  there exists a constant $\lambda= \lambda(\Phi)>1$, 
called the {\it dilatation of} $\Phi$, and there exists a pair of transverse measured foliations $\mathcal{F}^s$ and $\mathcal{F}^u$ such that 
 $$\Phi(\mathcal{F}^s)= \frac{1}{\lambda} \mathcal{F}^s \ \mbox{and}\  \Phi(\mathcal{F}^u)= \lambda \mathcal{F}^u.$$ 
 A mapping class  $\phi \in \mathcal{M}(\Sigma)$  is said to be {\it pseudo-Anosov} ({\it pA}) if  $\phi$ contains a pA homeomorphism. 
 We define the dilatation of a pA mapping class $\phi$, denoted by $\lambda(\phi)$, to be the dilatation of a pA homeomorphism of $\phi$.

Let $\mathcal{G}$ be a graph. 
We denote the set of vertices by $V(\mathcal{G})$  and denote the set of edges by $E(\mathcal{G})$. 
A continuous map $g: \mathcal{G} \rightarrow \mathcal{G}'$ from $\mathcal{G}$ into another graph $\mathcal{G}'$ is said to be a {\it graph map}. 
When $\mathcal{G}$ and $\mathcal{G}'$ are trees, a graph map $g: \mathcal{G} \rightarrow \mathcal{G}'$ is said to be a {\it tree map}. 
A graph map $g$ is called {\it Markov} if $g(V(\mathcal{G})) \subset V(\mathcal{G}')$ and 
for each point $x \in \mathcal{G}$ such that  $g(x) \notin V(\mathcal{G}')$, $g$ is locally injective at $x$ (that is $g$ has no $\lq$back track' at $x$). 
 In the rest of the paper we assume that all graph maps are Markov. 

 For a graph map $g$, we define the {\it transition matrix} 
 $M(g)=(m_{i,j})$ such that the $i^{\mathrm{th}}$ edge $e'_i$ or the same edge with opposite orientation  $(e'_i)^{-1}$  of $\mathcal{G}'$ appears 
 $m_{i,j}$-times in the edge path $g(e_j)$ for  the $j^{\mathrm{th}}$ edge $e_j$ of $\mathcal{G}$. 
If  $\mathcal{G}=\mathcal{G}'$, then $M(g)$ is a square matrix, and it makes sense to consider the 
 spectral radius, $\lambda(g) = \lambda(M(g))$, called the {\it growth rate} for $g$. 
 The topological entropy of $g$ is known to be equal to $\log \lambda(g)$.

In Section~\ref{subsection_PF}  we recall results regarding Perron-Frobenius matrices. 
In Section~\ref{subsection_traintrack} we quickly review a result from the train track theory which tells us that if 
a given mapping class $\phi$ induces a certain graph map, called {\it train track map}, whose transition matrix is Perron-Frobenius, 
then  $\phi$ is pA and $\lambda(\phi)$ equals the growth rate of the train track map. 
In Section~\ref{subsection_roots} 
we consider roots of a family of polynomials to study the dilatation of pA mapping classes and 
give some results regarding the asymptotic behavior of roots of this family.

\subsection{Perron-Frobenius theorem}
\label{subsection_PF}

Let $M=(m_{i,j}) $ and $N=(n_{i,j})$ be matrices with the same size. 
We shall write  $M \ge N$ (resp. $M >N$) whenever  $m_{i,j} \ge n_{i,j}$ (resp. $m_{i,j} > n_{i,j}$) for each $i,j$. 
We say that $M$ is {\it positive} (resp. {\it non-negative}) if $M >{\bf 0}$ (resp. $M \ge {\bf 0}$), where ${\bf 0}$ is the zero matrix.  

For a square and non-negative matrix $T$,  let $\lambda(T)$ be its spectral radius, that is the maximal absolute value of eigenvalues of $T$.  
We say that  $T$ is {\it irreducible} if for every pair of indices $i$ and $j$, there exists an integer $k=k_{i,j}>0$ such that the $(i,j)$ entry of 
$M^k$ is strictly positive. 
The matrix  $T$ is {\it primitive} if there exists an integer $k >0$ such that the matrix $T^k $ is positive.  
By definition, a primitive matrix is irreducible.  
A primitive matrix $T$ is {\it Perron-Frobenius}, abbreviated to PF,  if $T$ is an integral matrix. 
For $M \ge T$, if $T$ is irreducible then $M$ is also irreducible. 
The following theorem is commonly referred  to as  the Perron-Frobenius theorem.

\begin{thm} 
\cite{Sen}
\label{thm_PFtheorem}
Let $T$ be a primitive matrix. 
Then, there exists an eigenvalue $\lambda>0$ of $T$ such that 
\medskip
\\
(1) $\lambda$ has strictly positive left and right eigenvectors ${\bf \widehat{x}}$ and ${\bf y}$ respectively, and 
\medskip
\\
(2) $\lambda>|\lambda'|$ for any eigenvalue $\lambda' \ne \lambda$ of $T$. 
\end{thm}
\noindent
If $T$ is a PF matrix, the largest eigenvalue $\lambda$ in the sense of Theorem~\ref{thm_PFtheorem} 
is strictly greater than $1$, and it is called the {\it PF eigenvalue}. 
The corresponding positive eigenvector is called the {\it  PF eigenvector}.

The following will be  useful. 

\begin{lem} \cite[Theorem~1.6, Exercise~1.17]{Sen} 
\label{lem_subinvariance}
Let $T$ be a primitive matrix,  and let $s$ be a positive number. 
Suppose that   a non-zero vector ${\bf y} \ge {\bf 0}$ satisfies $T{\bf y} \ge s {\bf y}$. 
Then, 
\medskip
\\
(1) $\lambda \ge s$, where $\lambda$ is the largest eigenvalue  of $T$ in the sense of Theorem~\ref{thm_PFtheorem}, and 
\medskip
\\
(2) $s=\lambda$ if and only if $T{\bf y} = s {\bf y}$. 
\end{lem}

{\it Proof.} 
(1)   
Let $\widehat{{\bf x}}$ be a positive left eigenvector of $T$. 
Then, 
$$\widehat{{\bf x}} T {\bf y} =\lambda \widehat{{\bf x}} {\bf y} \ge s \widehat{{\bf x}} {\bf y}.$$
Hence we have $\lambda \ge s$. 

(2) ($\lq$Only if' part)  
Suppose that $s=\lambda$, and 
suppose that  $T{\bf y} \ge \lambda {\bf y}$ and $T{\bf y} \ne \lambda {\bf y}$. 
Premultiplying this inequality by  a positive left eigenvector $\widehat{{\bf x}}$ of $T$, 
we have 
$$\widehat{{\bf x}} T {\bf y}(=\lambda \widehat{{\bf x}} {\bf y})> \lambda \widehat{{\bf x}} {\bf y}.$$ 
Hence $\lambda>\lambda$, which is a contradiction.

($\lq$If' part)  
Suppose that $T{\bf y} = s {\bf y}$. 
Premultiplying this equality by  a positive left eigenvector $\widehat{{\bf x}}$ of $T$, 
we obtain $\lambda=s$. 
$\Box$ 
 \medskip

 For a non-negative $k \times k$ matrix $T$, one can associate a {\it directed graph} $G_T$ as follows. 
The graph $G_T$ has vertices numbered  $1, 2, \cdots, k$ and an edge from the $j^{\mathrm{th}}$ vertex  to the $i^{\mathrm{th}}$ vertex 
if and only if the $(i,j)$ entry $T_{i,j} \ne 0$. 
By the definition of $G_T$, one  easily verifies the following. 

\begin{lem} 
Let $T$ be a non-negative square matrix. 
\medskip
\\
(1) $T$ is irreducible if and only if for each $i,j$, there exists an integer $n_{i,j}>0$ such that 
the directed graph $G_T$ has an edge path of length $n_{i,j}$ from the $j^{\mathrm{th}}$ vertex to the $i^{\mathrm{th}}$ vertex. 
\medskip
\\
(2) $T$ is primitive if and only if there exists an  integer $n>0$ such that 
for each $i,j$, the directed graph $G_T$ has an edge path of length $n$   from the $j^{\mathrm{th}}$ vertex to the $i^{\mathrm{th}}$ vertex. 
\end{lem}

\subsection{Train track maps}
\label{subsection_traintrack}

A smooth branched $1$-manifold $\tau$ embedded in $D_n$ is a  {\it train track} if 
each component of $D_n \setminus \tau$ is either a non-punctured $k$-gon ($k \ge 3$), a once punctured $k$-gon ($k \ge1$) 
or an annulus such that  a boundary component of the annulus coincides with  the boundary of $D_n$ and the other component 
has at least  $1$ prong. 
A smooth map from a train track into itself is called a {\it train track map}.

Let $f: D_n \rightarrow D_n$ be a homeomorphism. 
A train track $\tau$ is {\it invariant} under $f$ if $f(\tau)$ can be collapsed smoothly onto $\tau$ in $D_n$. 
In this case $f$ induces a train track map $\widehat{f}: \tau \rightarrow \tau$.  
An edge of $\tau$ is called {\it infinitesimal} if there exists an integer $N>0$ such that $\widehat{f}^N(\tau)$ is a periodic edge under $\widehat{f}$. 
An edge of $\tau$ is called {\it real} if it is not infinitesimal. 
The transition matrix of $\widehat{f}$ is of the form:  
$$M(\widehat{f})= 
\left(\begin{array}{cc}
M_{\mathrm{real}}(\widehat{f}) & {\bf 0} 
\\
A & M_{\mathrm{inf}}(\widehat{f})
\end{array}\right),$$
where $M_{\mathrm{real}}(\widehat{f})$ (resp. $M_{\mathrm{inf}}(\widehat{f})$) is the transition matrix with respect to real (resp. infinitesimal) edges. 
The following is a consequence of \cite{BH}.

\begin{prop}
\label{prop_BH}
A mapping class  $ \phi \in \mathcal{M}(D_n)$ is pA if and only if 
there exists a homeomorphism $f: D_n \rightarrow D_n$ of $\phi$ and there exists a train track $\tau$ such that 
$\tau$ is invariant under $f$, and for the induced  train track map $\widehat{f}: \tau \rightarrow \tau$, the matrix $M_{\mathrm{real}}(\widehat{f})$ is PF. 
When  $\phi$ is a pA mapping class, we have $\lambda(\phi)=\lambda(M_{\mathrm{real}}(\widehat{f}))$. 
\end{prop}

 \subsection{Roots of polynomials}
 \label{subsection_roots}

For an integral polynomial $S(t)$, let $\lambda(S(t))$ be the maximal absolute value of roots of $S(t)$. 
For a monic integral polynomial $R(t)$, we set 
$$Q_{n,\pm}(t)= t^n  R(t) \pm S(t)$$ 
for each integer $n \ge 1$. 
The polynomial $R(t)$ (resp. $S(t)$)  is  called {\it dominant} (resp. {\it recessive}) for a family of polynomials $\{Q_{n,\pm}(t)\}_{n \ge 1} $. 
In case where $S(t)= R_*(t)$, we call $t^n  R(t) \pm R_*(t)$  the {\it Salem-Boyd polynomial} {\it associated to} $R(t)$. 
E.~Hironaka shows that such polynomials have several nice properties \cite[Section~3]{Hir}. 
The following lemma  shows that roots of $Q_{n,\pm}(t)$ lying outside the unit circle are determined by those of $R(t)$ asymptotically.

\begin{lem}
\label{lem_asymptotic-root1}
Suppose that $R(t)$ has a root outside the unit circle.  
Then, the roots of $Q_{n,\pm}(t)$ outside the unit circle converge to those of $R(t)$ counting multiplicity as $n$ goes to $\infty$.  
In particular, $ \lambda(R(t)) = \lim_{n \to \infty} \lambda(Q_{n,\pm}(t))$. 
\end{lem}

\noindent
The proof can be found in \cite{Hir}. 
We recall a proof here for completeness. 
\medskip

{\it Proof.} 
Consider the rational function 
$$ \frac{Q_{n,\pm}(t)}{t^n}= R(t) \pm \frac{S(t)}{t^n}.$$
Let $\theta$ be a root of $R(t)$ with multiplicity $m$ outside the unit circle. 
Let $D_{\theta}$ be any small disk centered at $\theta$ that is strictly outside of the unit circle 
and that contains no roots of $R(t)$ other than $\theta$. 
Then, $|R(t)|$ has a lower bound on the boundary $\partial D_{\theta}$ by compactness. 
Hence there exists a number $n_{\theta}>0$ depending on $\theta$ such that $|R(t)|> |\frac{S(t)}{t^n}|$ on $\partial D_{\theta}$ for any $n > n_{\theta}$. 
By the Rouch\'{e}'s theorem, it follows that $R(t)$ and $R(t) \pm \frac{S(t)}{t^n}$ (hence $R(t)$ and $Q_{n,\pm}(t)$)  have the same $m$ roots in $D_{\theta}$. 
Since $D_{\theta}$ can be made arbitrarily small and there exist only finitely many roots of $R(t)$, 
the proof of Lemma~\ref{lem_asymptotic-root1} is complete. 
$\Box$

\begin{lem}
\label{lem_asymptotic-root2}
Suppose that $R(t)$ has no roots outside the unit circle, and suppose that 
$Q_{n,\pm}(t)$ has  a real root $\mu_n$ greater than $1$ for sufficiently large $n$. 
Then, $\lim_{n \to \infty} \mu_n= 1$. 
\end{lem}

{\it Proof.} 
For any $\varepsilon>0$, let $D_{\varepsilon}$ be the disk of radius $1+ \varepsilon$ around the origin in the complex plane. 
Then, for any sufficiently large $n$, we have $|R(t)| > |\frac{S(t)}{t^n}|$  for all $t$ on $\partial D_{\varepsilon}$. 
Moreover, $R(t)$ and $\pm \frac{S(t)}{t^n}$ are holomorphic on the complement of $D_{\varepsilon}$ in the Riemann sphere. 
By Rouch\'e's theorem, $R(t)$ and $R(t) \pm \frac{S(t)}{t^n}$ 
(hence $R(t)$ and $Q_{n,\pm}(t)$) have no roots outside $ D_{\varepsilon}$. 
Hence $\mu_n$ converges to $1$ as $n$ goes to $\infty$. 
$\Box$

 \section{Combined tree maps}
 \label{section_combined}

For an $n \times n$ matrix $M$, let $M(t)$ be the  characteristic polynomial $|tI-M|$ of $M$, where $ I = I_n$ is the $n \times n$ identity matrix. 
Let $M_*(t)$ be the reciprocal polynomial of $M(t)$. 
Then, 
$$M_*(t)= t^n \Bigl|\frac{1}{t} I - M \Bigr|=   |I - tM|,$$ 
that is $M_*(t)$ equals the determinant of the matrix $I - tM$.

This section introduces  {\it combined tree maps}. 
Given two trees we combine these trees with another tree of star type having the valence $n+1$ vertex 
and define a new tree, say $\mathcal{Q}_n$. 
When two tree maps on  $\mathcal{Q}_n$ satisfy certain conditions ({\bf L1},{\bf L2},{\bf L3} and {\bf R1},{\bf R2},{\bf R3}),  
we can define the combined tree map $q_n$ on $\mathcal{Q}_n$ and obtain a family of tree maps $\{q_n: \mathcal{Q}_n \rightarrow \mathcal{Q}_n\}_{n \ge 1}$.  
In Section~\ref{subsection_growth} we give a sufficient condition that guarantees $M(q_n)$ is PF. 
In Section~\ref{subsection_asymptotic} we consider combined tree maps in a particular setting. 
Then, we give a formula for $M(q_n)(t)$ and ${M(q_n)}_*(t)$ and analyze the asymptotic behavior of the growth rate for $q_n$. 
This analysis  will be  applied to train track maps in Section~\ref{section_proof}.

 \subsection{Transition matrices and growth rate} 
 \label{subsection_growth}

\begin{figure}[htbp]
\begin{center}
\includegraphics[width=2.5in]{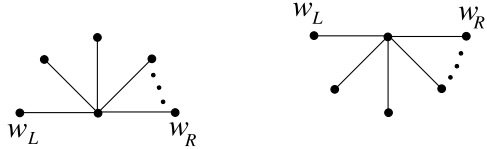}
\caption{trees (left) $\mathcal{G}_{n,+}$ and (right) $\mathcal{G}_{n,-}$ having one vertex of valence $n+1$. }
\label{fig_star}
\end{center}
\end{figure}

 We assume that all trees are embedded in the disk $D$. 
 By the trivial tree $\mathcal{T}_0$, we mean  the tree with only one vertex. 
 Let $\mathcal{G}_{n,+}$ and $\mathcal{G}_{n,-}$ be trees of star type  as in Figure~\ref{fig_star}, having one vertex of valence $n+1$.

\begin{figure}[htbp]
\begin{center}
\includegraphics[width=4.5in]{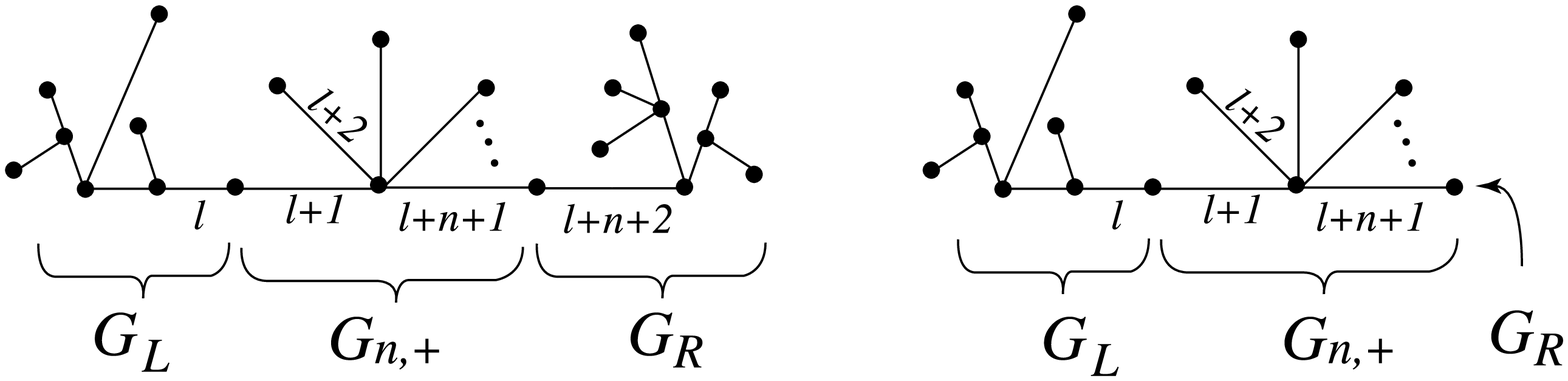}
\caption{combined trees  $\mathcal{Q}_{n,+}$: 
(left) general case, (right) case where $\mathcal{G}_R$ is a trivial tree.}
\label{fig_tree_Q}
\end{center}
\end{figure}

 Let $\mathcal{G}_L$ (resp. $\mathcal{G}_R$) be  a tree (possibly a trivial tree)  
 with a valence $1$ vertex, say $v_L$ (resp. $v_R$). 
 Let  $w_L$ and $w_R$ be vertices of $\mathcal{G}_{n,+}$ as in Figure~\ref{fig_star}, and glue $\mathcal{G}_L$, $\mathcal{G}_{n,+}$ and 
 $\mathcal{G}_R$ together so that for $S \in \{L,R\}$, $v_S $ and $w_S$ become one vertex (Figure~\ref{fig_tree_Q}).  
The resulting tree $\mathcal{Q}_{n,+}$ is called the {\it combined tree}, obtained from the triple ($\mathcal{G}_L, \mathcal{G}_{n,+}, \mathcal{G}_R$). 
We define the combined tree $\mathcal{Q}_{n,-}$, obtained from the triple   ($\mathcal{G}_L, \mathcal{G}_{n,-}, \mathcal{G}_R$) in the same manner.

\begin{figure}[htbp]
\begin{center}
\includegraphics[width=4in]{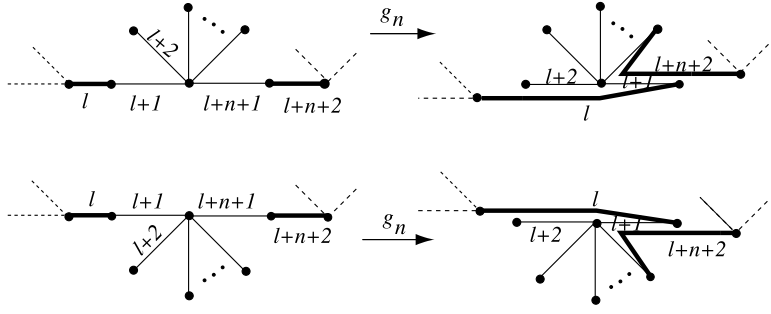}
\caption{(top) $g_n$ rotates $\mathcal{G}_{n,+}$, (bottom) $g_n$ rotates $\mathcal{G}_{n,-}$. 
The edges $e_{\ell}$ and $e_{\ell+n+2}$ and their images are drawn in bold.} 
\label{fig_g_n}
\end{center}
\end{figure}

Before we define combined tree maps on $\mathcal{Q}_{n,+/-}$, 
we label  the  edges of $\mathcal{Q}_{n,+/-}$. 
Let  $\ell$ be  the number of edges of $\mathcal{G}_L$, and let  $r$ be  the number of edges of $\mathcal{G}_R$ plus $1$. 
Note that the number of edges of $\mathcal{Q}_{n,+/-}$ is $\ell + n+ r$. 
\medskip
\\
$\bullet$
The edges of  $\mathcal{G}_{n,+/-}$ are numbered $\ell+1$ to $\ell+n+1$ in the clockwise/counterclockwise direction as in 
Figure~\ref{fig_g_n}.  
\medskip
\\
$\bullet$ 
The edge of $\mathcal{G}_L$ sharing a vertex with the $(\ell+1)^{\mathrm{st}}$ edge is numbered $\ell$ 
and the remaining  edges of $\mathcal{G}_L$ are numbered $1$ to $\ell-1$ arbitrarily. 
\medskip
\\
$\bullet$ 
The edge of $\mathcal{G}_R$  sharing a vertex with the $(\ell+n+1)^{\mathrm{st}}$ edge is numbered  
$\ell+n+2$ and the remaining  edges of $\mathcal{G}_R$ are numbered $\ell+n+3$ to $ \ell+n+r$ arbitrarily. 
\medskip
\\
The edge numbered $i$ is denoted by $e_i$.

Now we take a tree map  $g_L: \mathcal{Q}_{n,+/-} \rightarrow \mathcal{Q}_{n,+/-}$ satisfying the following conditions. 
\medskip
\\
{\bf L1}
The map $g_L$ restricted to the set of vertices  of $E(\mathcal{Q}_{n,+/-}) \setminus(E(\mathcal{G}_L) \cup \{e_{\ell+1}\})$ is the identity. 
\medskip
\\
{\bf L2}
$g_L(\mathcal{G}_L) \subset \mathcal{G}_L$. 
\medskip
\\
{\bf L3}
The edge path $g_L(e_{\ell+1})$ passes through  $e_{\ell+1}$  only once and passes through $e_{\ell}$. 
\medskip

Next, we take a tree map $g_R: \mathcal{Q}_{n,+/-} \rightarrow \mathcal{Q}_{n,+/-}$ satisfying  the following conditions. 
\medskip
\\
{\bf R1}
The map $g_R$ restricted to the set of vertices of $E(\mathcal{Q}_{n,+/-}) \setminus(E(\mathcal{G}_R) \cup \{e_{\ell+n+1}\})$ is the identity. 
 \medskip
 \\
{\bf R2} 
$g_R(\mathcal{G}_R) \subset \mathcal{G}_R$. 
\medskip
\\
{\bf R3}
The edge path $g_R(e_{\ell+n+1})$ passes through $e_{\ell+n+1}$ only once and passes through $e_{\ell +n+2}$.  
\medskip

Finally, we define the  tree map   $g_n: \mathcal{Q}_{n,+/-} \rightarrow \mathcal{Q}_{n,+/-}$ satisfying  the following conditions. 
\medskip
\\
{\bf n1}
The map $g_n$ restricted to the set of vertices of 
$E(\mathcal{Q}_{n,+/-}) \setminus(E(\mathcal{G}_{n,+/-}) \cup \{e_{\ell} , e_{\ell+n+2}\})$ is the identity. 
 \medskip
 \\
 {\bf n2}
 $g_n$ rotates the subtree $\mathcal{G}_{n,+/-}$ as in Figure~\ref{fig_g_n}. 
 \medskip
 \\
{\bf n3}
 The image of  each  $e \in \{e_{\ell} , e_{\ell+n+2}\}$ is as in Figure~\ref{fig_g_n}.  
The length of the edge path $g_n(e)$ is $3$.
\medskip

The composition 
$$q_n= g_R  g_{n}  g_L : \mathcal{Q}_{n,+/-} \rightarrow \mathcal{Q}_{n,+/-}$$ 
is called the  {\it combined tree map}, obtained from the triple $(g_L, g_n, g_R)$. 
It makes sense to consider the transition matrices $M(g_S)$ of 
$g_S|_{\mathcal{G}_S}: \mathcal{G}_S  \rightarrow \mathcal{G}_S$, $S \in \{L,R\}$ and 
$M(g_n)$ of $g_n|_{\mathcal{G}_{n,+/-}}: \mathcal{G}_{n,+/-} \rightarrow \mathcal{G}_{n,+/-}$. 
The transition matrix $M(q_n)$ has the following form: 
\begin{equation}
\label{equation_matrix}
 M(q_{n})= 
 \left(\begin{array}{ccc}
 M_L   & A & {\bf 0}\\
 B &M_n & C\\
 D & E & M_R
 \end{array}\right),
 \mbox{where} \ 
 M_n=
 \begin{scriptsize}
 \left(\begin{array}{cccc} * & 1 &  &  \\
  &  & \ddots &  \\
   &  &  & 1 \\
  * &  &  & 
   \end{array}\right)
   \end{scriptsize}
\end{equation}
(each empty space in $M_n$ represents the number $0$),  
and the block matrices satisfy $M_L \ge M(g_L)$, $M_n \ge M(g_n)$ and $M_R \ge M(g_R)$. 
(In fact $M_L = M(g_L)$, although  we will not be using  this fact.) 

Throughout this subsection, we assume that the trees $\mathcal{G}_L$ and $\mathcal{G}_R$ are not trivial.  
It is straightforward to see the following from the defining conditions of $g_L$, $g_n$ and $g_R$.

\begin{lem}
\label{lem_matrix}
Let $m_{i,j}$ be the $(i,j)$ entry of $M(q_n)$. We have 
\\
(1) 
 $m_{\ell+n,\ell+n+2}=1$ and $ m_{\ell+n+1,\ell+n+2}=1$, and 
\medskip
\\
(2) 
$m_{\ell,\ell+1}>0$, $m_{\ell+1,\ell+1}>0$ and $m_{\ell+n+1,\ell+1}>1$. 
Moreover, $m_{\ell+1,j}= m_{\ell+n+1,j}$ for each $j$ with $1 \le j \le \ell$ and 
$m_{\ell+1,j_0}>0$ for some $1 \le j_0 \le \ell$, \mbox{and}
\medskip
\\
(3) $m_{\ell+n+2,\ell+1}>0$. 
\end{lem}

An important feature  is that the growth rate of $q_n$ is always greater than $1$ if $M(g_L)$ and $M(g_R)$ are irreducible, which will be shown  
in Proposition~\ref{prop_PF-matrix}. 
We first show that $M(q_n)$ is irreducible in this case. 
Notice that $M(g_n)$ is always irreducible, and since $M_n \ge M(q_n)$ so $M_n$ must be irreducible as well.

\begin{lem}
\label{lem_irreducible}
Let $q_n= g_R  g_{n}  g_L: \mathcal{Q}_{n,+/-} \rightarrow \mathcal{Q}_{n,+/-}$ be the combined tree map. 
Assume that  both $M(g_L)$ and $M(g_R)$ are irreducible. Then, $M(q_n)$ is irreducible.
\end{lem}

{\it Proof.} 
Note that $M_L$ $M_R$ and $M_n$ are irreducible. 
Let $G_{q_n}$ be the directed graph of $M(q_n)$. 
We identify vertices of $G_{q_n}$ with edges of $\mathcal{Q}_{n,+/-}$. 
Let $V_L$ (resp. $V_R$, $V_n$) be the set of vertices  of $G_{q_n}$ coming from the set of edges of the subtree $\mathcal{G}_L$ 
(resp. $\mathcal{G}_R$, $\mathcal{G}_{n,+/-}$) of $\mathcal{Q}_{n,+/-}$. 
Lemma~\ref{lem_matrix}(2) shows that there exists an edge connecting the set $V_L$ to the set $V_n$, and 
there exists an edge connecting the set $V_n$ to the set $V_L$. 
This is also  true between $V_n$ and $V_R$ by  Lemma~\ref{lem_matrix}(1,3). 
Thus, one can find   an edge path between any two vertices of $G_{q_n}$. 
$\Box$

\begin{prop}
\label{prop_PF-matrix}
Under the assumptions of Lemma~\ref{lem_irreducible}, $M(q_{n})$ is  PF. 
\end{prop}

{\it Proof.} 
Lemma~\ref{lem_matrix}(2) says that the directed graph $G_{q_n}$ has an edge from the vertex $v_{\ell+1}$ to itself, and we denote such edge by $e$. 
Since $M(q_n)$ is irreducible, for any vertex $v$ of $G_{q_n}$ 
there exists an edge path $E= e_1 e_2 \cdots e_{n(v)}$ from $v_{\ell+1}$ to $v$. 
Thus,  for any $n \ge n(v)$ we have an edge path $e \cdots e  E$ of length $n$ from $v_{\ell+1}$ to $v$. 
Since the number of vertices is finite,  there exists an integer $N>0$ such that for any vertex $w$ of $G_{q_n}$ and any integer $n \ge N$ 
we have an edge path of length $n$ from $v_{\ell+1}$ to $w$. 
Since there exists an edge path from any vertex $x$ of $G_{q_n}$ to $v_{\ell+1}$, 
we can find a sufficiently large integer $N'$ such that for any pair of vertices $x$ and  $w$ there exists an edge path of length $N'$ from $x$ to $w$. 
Thus, $M(q_n)$ is PF. 
$\Box$ 
\medskip

The following property is crucial in proving Proposition~\ref{prop_decreasing-example} and Theorem~\ref{thm_lim-example}.

\begin{prop}
\label{prop_decreasing}
Under the assumptions of Lemma~\ref{lem_irreducible}, we have $\lambda(M(q_{n})) >  \lambda(M(q_{n+1}))>1$. 
\end{prop}
 
 {\it Proof.} 
To compare  $M(q_{n+1})$ with $M(q_n)$ we introduce a new labeling of edges of $\mathcal{Q}_{n+1,+/-}$. 
The trees $\mathcal{G}_L$ and $\mathcal{G}_R$ are the common subtrees for both trees $\mathcal{Q}_{n,+/-}$ and $\mathcal{Q}_{n+1,+/-}$. 
Edges of  the subtrees $\mathcal{G}_L$ and $\mathcal{G}_R$ of $\mathcal{Q}_{n+1,+/-}$ are numbered in the same manner as those of $\mathcal{Q}_{n,+/-}$, and 
 edges  of $\mathcal{G}_{n+1}$ are numbered  
 $$\ell+1, \ell+n+r+1, \ell+2, \ell+3, \cdots, \ell+n+1$$
 in the clockwise/counterclockwise direction. 
Here the edge sharing a  vertex with the $\ell^{\mathrm{th}}$ edge is numbered $\ell+1$. 
 
 Let $M(q_n)= (m_{i,j})_{1 \le i,j \le \ell+n+r}$ be the matrix given in (\ref{equation_matrix}). 
Then, $M(q_{n+1} )= (m'_{i,j})_{1 \le i,j \le \ell+n+r+1}$ with new labeling has the following form: 
\medskip

 \begin{center}
 $M(q_{n+1})=$ 
 \begin{scriptsize}
 $
 \left(
 \begin{array}{cc|ccccc|ccc|c}
  M_L      & & && A  & & & & {\bf 0}   & & \\ 
                                      & & & &                   & & & &                            & & \\ \hline 
                                  & &* & &                    & & & &                             & &  1\\ 
                                  & & & & 1                   & & & &                            & & \\ 
  B      & & & &                     &\ddots& & & C         & & \\ 
                                 & & & &                      &           &1 & &                           & & \\ 
                                   & &* & &                     & & & &                            & & \\ \hline  
                                    & & & &                    & & & &                            & & \\ 
   D            & & & & E& & &  & M_R     & & \\ 
                                    & & & &                     & & & &                            & & \\ \hline   
                                     & & &1 &                    & & & &                              & & \\ 
 \end{array}
 \right).
 $ 
\end{scriptsize}
\medskip
\end{center}
Put $s= \lambda(M(q_{n+1}))>1$ and let ${\bf y}= \ ^t(y_1, \cdots, y_{\ell+n+r+1})$ be the PF eigenvector for $M(q_{n+1})$. 
Then, 
\begin{equation}
\label{equation_PF-equality} 
\sum_{j=1}^{\ell+n+r+1} m'_{i,j}y_j= sy_j \hspace{2mm} \mbox{for\ }i\ \mbox{with}\ 1 \le i \le \ell+n+r+1.
\end{equation}
For $i=\ell+1$ and $i= \ell+n+r+1$ of (\ref{equation_PF-equality}) we have 
\begin{eqnarray*}
\sum_{j=1}^{\ell+n+r+1} m'_{\ell+1,j}y_j &=& \sum_{j=1}^{\ell+1} m_{\ell+1,j} y_j+ y_{\ell+n+r+1} = sy_{\ell+1}\ \mbox{and}
\\
\label{equation_PF-equality2}
y_{\ell+2} &=& sy_{\ell+n+r+1}. 
\end{eqnarray*}
These two equalities together with $s>1$ yield 
\begin{equation}
\label{equation_PF-equality3}
\sum_{j=1}^{\ell+1} m_{\ell+1,j} y_j+ y_{\ell+2} > sy_{\ell+1}. 
\end{equation}
The equalities~(\ref{equation_PF-equality}) for all $i \ne \ell+1,\ell+n+r+1$ together with 
the inequality (\ref{equation_PF-equality3}) imply 
$$M(q_n) \hat{{\bf y}} \ge s\hat{{\bf y}},\ \mbox{where}\ \hat{{\bf y}}= \ ^t(y_1, \cdots, y_{\ell+n+r}).$$
By Lemma~\ref{lem_subinvariance}(1), we have $\lambda(M(q_{n})) \ge  s= \lambda(M(q_{n+1}))$.  
By Lemma~\ref{lem_subinvariance}(2) together with (\ref{equation_PF-equality3}), we have $\lambda(M(q_{n})) >  s$. 
$\Box$

 \subsection{Asymptotic behavior of growth rate}
 \label{subsection_asymptotic}

In this section we concentrate on  the combined tree obtained from the triple $(\mathcal{G}_L, \mathcal{G}_{n,+/-}, \mathcal{T}_0)$. 
We assume that  $g_L(\mathcal{G}_L)= \mathcal{G}_L$ and study the combined tree map $q_n= g_n g_L$.

\begin{figure}[htbp]
\begin{center}
\includegraphics[width=4in]{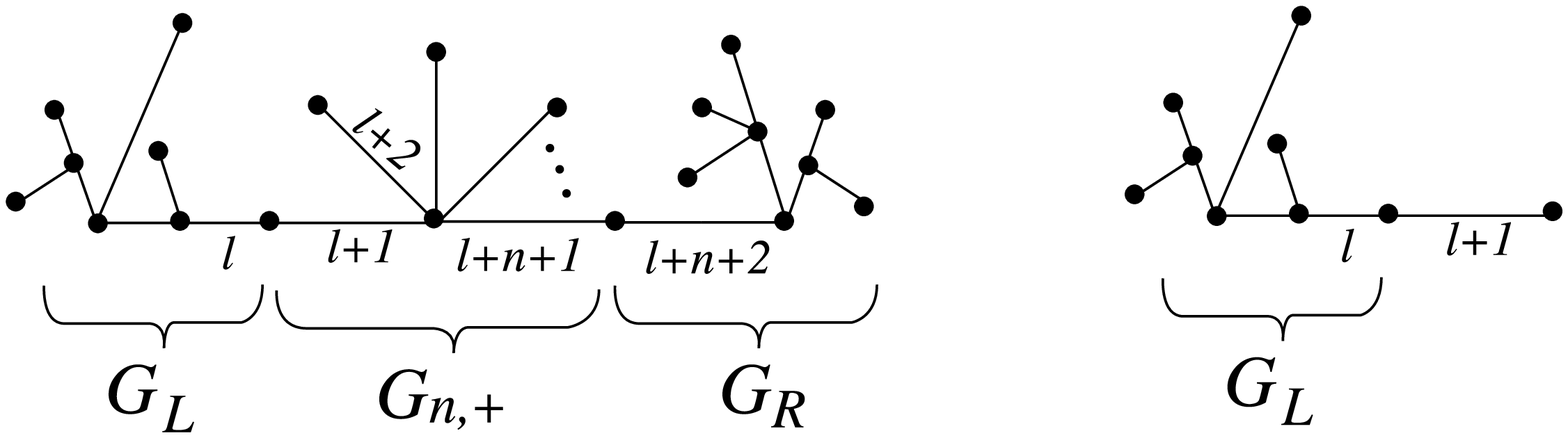}
\caption{(left) ${\mathcal Q}_{n,+}$, (right) its subtree ${\mathcal R}$} 
\label{fig_subtree_R}
\end{center}
\end{figure}

Let ${\mathcal R}$ be the subtree  of ${\mathcal Q}_{n,+/-}$ such that ${\mathcal R}$ is obtained from the subtree 
${\mathcal G}_L$ together with the $(\ell+1)^{\mathrm{st}}$ edge. (Hence $E({\mathcal R})=\{e_1, e_2, \cdots, e_{\ell+1}\}$.) 
For an example of ${\mathcal R}$, see Figure~\ref{fig_subtree_R}.

By the assumption $g_L(\mathcal{G}_L)= \mathcal{G}_L$, we have $q_n(\mathcal{R}) \supset \mathcal{R}$ 
and hence the following  tree map  $\overline{r}: \mathcal{R} \rightarrow \mathcal{R}$ is well defined:  
for each $e \in E({\mathcal R})$,  the edge path $\overline{r}(e)$ is given by  the edge path $q_n(e)$  by eliminating edges which do not belong to $E(\mathcal{R})$. 
The tree map $\overline{r}$ does not depend on the choice of $n$. 
The transition matrix $M(\overline{r})$ is given by  the upper-left $(\ell+1) \times (\ell+1)$ submatrix of $M(q_n)$.  
We call  ${\mathcal R}$ the {\it dominant tree} and $\overline{r}$ the {\it dominant tree map} for a  family of combined tree maps $\{q_{n}\}_{n  \ge 1}$.

We now define a polynomial  $S(t)$ (resp. $U(t)$) as follows: 
Consider the matrix $tI - M(q_n)$ (resp. $I - t M(q_n)$) and replace the $(\ell+1)^{\mathrm{st}}$ row by the last row. 
Take the upper-left  $(\ell+1) \times (\ell+1)$ submatrix of the resulting matrix, denoted by $S$ (resp. $U$), and then  
 $S(t)$ (resp. $U(t)$) is defined equal to the determinant of $S$ (resp. $U$). 
It is not hard to see  that the matrices $S$  and $U$ do not depend on  $n$.

The following statement, which will be crucial later, tells us that  
$M(\overline{r})(t)$ is the dominant polynomial and $S(t)$ is the recessive polynomial for a family of polynomials $\{M(q_n)(t)\}_{n \ge 1}$. 

\begin{prop}
\label{prop_recipe}
We have 
\medskip
\\
(1)  $M(q_{n})(t)= t^n M(\overline{r})(t)+ S(t)$, and 
 \medskip
 \\
 (2) $M(q_n)_*(t) = t^nU(t) + M(\overline{r})_*(t)$. 
 \end{prop}
 
 {\it Proof.} 
The transition matrix $M(q_{n})= (m_{i,j})$ is of the form

$$
M(q_{n})= 
\begin{scriptsize}
\left(\begin{array}{ccc|ccc}  
&   & & & & 
\\
& M(\overline{r})& & & & 
\\
& & &1 & & 
\\
\hline 
& & & & \ddots& 
\\
& & & & & 1
\\
m_{\ell+n+1,1}& \cdots &m_{\ell+n+1,\ell+1} & & &
\end{array}\right), 
\end{scriptsize}\ 
$$
and it is easy to see that $m_{\ell+1,j}= m_{\ell+n+1,j}$ for $1 \le j \le \ell$. 
For the proof of (1) (resp. (2)), apply the determinant expansion with respect to  the last row of $tI-M(q_n)$ (resp. $I-t M(q_n))$. 
 $\Box$

\begin{prop}
\label{prop_PF}
Suppose that $M(g_L)$ is irreducible. 
Then, we have 
\medskip
\\
(1) $M(q_{n})$ is PF  for each $n$ and $\lambda(M(q_{n})) > \lambda(M(q_{n+1}))$, and 
\medskip
\\
(2)  $\lambda(M(\overline{r}))= \lim_{n \to \infty} \lambda(M(q_{n}))$. 
 \end{prop}
 
 {\it Proof.} 
 (1) The proof  is parallel to the proofs of Propositions~\ref{prop_PF-matrix} and \ref{prop_decreasing}. 
 
 (2) Apply Lemma~\ref{lem_asymptotic-root1} with Proposition~\ref{prop_recipe}(1). 
 $\Box$

\section{Proof}
\label{section_proof}

This section is devoted to proving Proposition~\ref{prop_decreasing-example} and Theorems~\ref{thm_recursive-example}, \ref{thm_lim-example}, \ref{thm_dil-vol}.

\begin{prop}
\label{prop_pseudo-Anosov}
The braids $\beta_{(m_1, \cdots, m_{k+1})}$ are pA. 
\end{prop}

\begin{figure}[htbp]
\begin{center}
\includegraphics[width=3in]{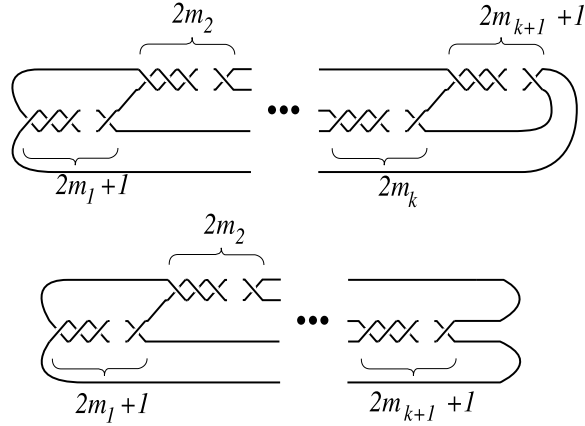}
\caption{$\overline{\beta_{(m_1, \cdots, m_{k+1})}}$: (top) $k$ odd, (bottom) even.}
\label{fig_general_2bridge}
\end{center}
\end{figure}

\noindent
By a result of W.~Menasco's \cite[Corollary~2]{Men}, if $L$ is a non-split prime alternating link which is not a torus link, then $S^3 \setminus L$ has a complete hyperbolic structure of finite volume. 
Since $\overline{\beta_{(m_1, \cdots, m_{k+1})}}$ is a $2$ bridge link as depicted in Figure~\ref{fig_general_2bridge}, his result tells us that 
$\beta_{(m_1, \cdots, m_{k+1})}$ is pA. 
Here we will show Proposition~\ref{prop_pseudo-Anosov} by using Proposition~\ref{prop_BH}. 
As a result, we will find the polynomial whose largest root equals the dilatation of $\beta_{(m_1, \cdots, m_{k+1})}$. 
\medskip

{\it Proof of Proposition~\ref{prop_pseudo-Anosov}.} 
To begin with,  we define a tree $\mathcal{Q}_{(m_1, \cdots, m_{k+1})}$ and a tree map $q_{(m_1, \cdots, m_{k+1})}$ on the tree 
$\mathcal{Q}_{(m_1, \cdots, m_{k+1})}$ inductively. 

\begin{figure}[htbp]
\begin{center}
\includegraphics[width=4in]{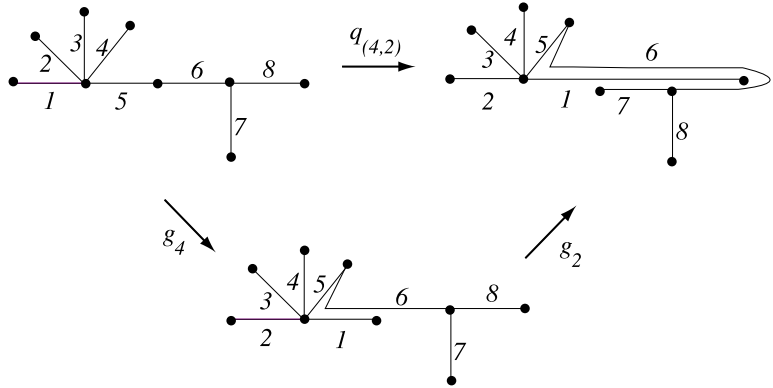}
\caption{case $(m_1, m_2) = (4,2)$. } 
\label{fig_combined_treemap}
\end{center}
\end{figure}

For $k = 1$ 
let $\mathcal{Q}_{(m_1,m_2)}$ be the combined tree obtained from the triple $(\mathcal{G}_{m_{1,+}}, \mathcal{G}_{m_2,-},\mathcal{T}_0)$. 
Take the tree maps $g_{m_1}$ and $g_{m_2}$ with conditions {\bf n1},{\bf n2},{\bf n3} and let us define $q_{(m_1,m_2)}$ as the combined tree map 
$$q_{(m_1,m_2)}= g_{m_2} g_{m_1}: \mathcal{Q}_{(m_1,m_2)} \rightarrow \mathcal{Q}_{(m_1,m_2)}\hspace{2mm} (\mbox{Figure}~\ref{fig_combined_treemap}).  $$

Next, suppose that  these are defined up to $k$. 
Let $\mathcal{Q}_{(m_1, \cdots, m_{k+1})}$ be the combined tree obtained from the triple 
$(\mathcal{Q}_{(m_1, \cdots, m_k)}, \mathcal{G}_{m_{k+1},+/-}, \mathcal{T}_0)$ in case $k+1$ odd/even.  
We extend $q_{(m_1, \cdots, m_k)}: \mathcal{Q}_{(m_1, \cdots, m_k)} \rightarrow \mathcal{Q}_{(m_1, \cdots, m_k)}$ 
 to a tree map 
 $$\widehat{q}: \mathcal{Q}_{(m_1, \cdots, m_{k+1})} \rightarrow \mathcal{Q}_{(m_1, \cdots, m_{k+1})}$$ 
 satisfying {\bf L1},{\bf L2},{\bf L3} 
so that  for the edge $e$ of $\mathcal{G}_{m_{k+1},+/-}$ sharing  a vertex with the edge of $\mathcal{Q}_{(m_1, \cdots, m_k)}$,  
the length of the edge path $\widehat{q}(e)$ is $3$. 
Let  us define  the combined tree map
$$q_{(m_1, \cdots, m_{k+1})}= g_{m_{k+1}}  \widehat{q}: \mathcal{Q}_{(m_1, \cdots, m_{k+1})} \rightarrow \mathcal{Q}_{(m_1, \cdots, m_{k+1})} .$$
By Proposition~\ref{prop_PF}(1), the transition matrix $M(q_{(m_1, \cdots, m_{k+1})})$ is PF.

\begin{figure}[htbp]
\begin{center}
\includegraphics[width=4in]{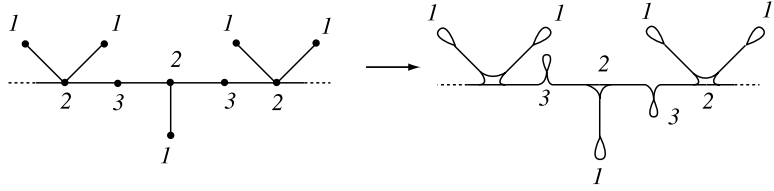}
\caption{numbers $1,2$ and $3$ correspond to deformations $1,2$ and $3$.} 
\label{fig_train-track}
\end{center}
\end{figure}

We now deform $\mathcal{Q}=\mathcal{Q}_{(m_1, \cdots, m_{k+1})}$ into a train track $\tau_{(m_1, \cdots, m_{k+1})}$ (as in Figure~\ref{fig_train-track}): 
\medskip
\\
1. 
Puncture a disk near each valence $1$ vertex of $\mathcal{Q}$, and connect a $1$-gon at the  vertex which contains the puncture. 
\medskip
\\
2. 
Deform a neighborhood of  a valence $m_i+1$ vertex of  the subtree $\mathcal{G}_{m_i,+/-}$ of $\mathcal{Q}$ into an $(m_i+1)$-gon for each $i$.  
 \medskip
 \\
3.  
For each $i+1$ odd/even, 
puncture above/below the vertex of valence $2$ which connects the two subtrees  $\mathcal{G}_{m_i,+/-}$ and $\mathcal{G}_{m_{i+1},-/+}$. 
Deform a neighborhood of the vertex and  connect a $1$-gon which contains the puncture.

\begin{figure}[htbp]
\begin{center}
\includegraphics[width=4in]{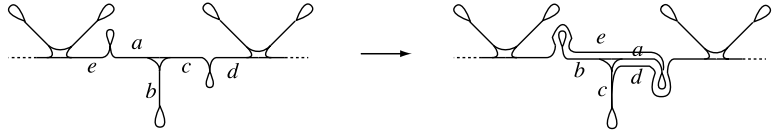}
\caption{$\widehat{q}$ rotates a part of the  train track smoothly.} 
\label{fig_train-track-map}
\end{center}
\end{figure}

Then, $q=q_{(m_1, \cdots, m_{k+1})}$ induces the graph map $\widehat{q}$ 
on $\tau= \tau_{(m_1, \cdots, m_{k+1})}$ into itself. 
Since $\widehat{q}$ rotates a part of the train track smoothly (Figure~\ref{fig_train-track-map}), it turns out  that $\widehat{q}$ is a smooth map. 
It is easy to show the existence of a representative homeomorphism $f$ of $\phi= \Gamma(\beta_{(m_1, \cdots, m_{k+1})})$ 
such that $\tau$ is invariant under $f$ and such  that 
$\widehat{q}$ is the  train track map induced by $f$. 
The transition matrix of $\widehat{q}$ with respect to real edges of $\tau$ is exactly equal to  the PF matrix $M(q_{(m_1, \cdots, m_{k+1})})$. 
By Proposition~\ref{prop_BH}, the braid $\beta_{(m_1, \cdots, m_{k+1})}$ is pA. 
This completes the proof of Proposition~\ref{prop_pseudo-Anosov}. 
$\Box$

\begin{figure}[htbp]
\begin{center}
\includegraphics[width=4in]{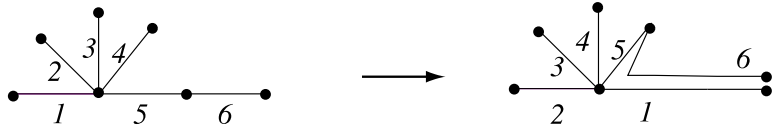}
\caption{dominant tree map $\overline{r}: \mathcal{R} \rightarrow \mathcal{R}$ for $\{q_{(4,m)}\}_{m \ge 1}$.} 
\label{fig_dominant_treemap}
\end{center}
\end{figure}

\begin{ex}
Let us express the formula to compute the dilatation of the braids $\beta_{(4,m)}$ for $m \ge 1$.  
We know that the dilatation of the braid $\beta_{(4,m)}$ is the largest root of the polynomial $M(q_{(4,m)})(t)$ for the combined tree map 
$q_{(4,m)}$. 
By using $q_{(4,2)}$ shown in Figure~\ref{fig_combined_treemap}, the transition matrix for $q_{(4,2)}$ is 
$$M(q_{(4,2)})= 
\begin{scriptsize}
\left(\begin{array}{cccccccc}
0 & 1 & 0 & 0 & 0 & 0 & 0 & 0 \\
0 & 0 & 1 & 0 & 0 & 0 & 0 & 0 \\
0 & 0 & 0 & 1 & 0 & 0 & 0 & 0 \\
0 & 0 & 0 & 0 & 1 & 1 & 0 & 0 \\
1 & 0 & 0 & 0 & 0 & 1 & 0 & 0 \\
1 & 0 & 0 & 0 & 0 & 1 & 1 & 0 \\
0 & 0 & 0 & 0 & 0 & 0 & 0 & 1 \\
1 & 0 & 0 & 0 & 0 & 2 & 0 & 0
\end{array}\right). 
\end{scriptsize}$$
The transition matrix $M(\overline{r})$ of the dominant tree map $\overline{r}$ (Figure~\ref{fig_dominant_treemap}) for a family of combined tree maps 
$\{q_{(4,m)}\}_{m \ge 1}$ is the upper-left $6 \times 6$ matrix. 
Hence the dominant polynomial $M(\overline{r})(t)$ for $\{q_{(4,m)}\}_{m \ge 1}$ equals $ t^6-t^5-2t $ with the largest root $\approx 1.45109$. 
In this case the recessive polynomial $S(t)$ is 
$$S(t)= 
\begin{scriptsize}
\left|\begin{array}{cccccc}
t & -1 & 0 & 0 & 0 & 0 \\
0 & t & -1 & 0 & 0 & 0 \\
0 & 0 & t & -1 & 0 & 0 \\
0 & 0 & 0 & t & -1 & -1 \\
-1 & 0 & 0 & 0 & t & -1 \\
-1 & 0 & 0 & 0 & 0 & -2
\end{array}\right| 
\end{scriptsize}= -2t^5-t+1.$$
By Proposition~\ref{prop_recipe}(1), the dilatation of $\beta_{(4,m)}$ is the largest root of 
$$M(q_{(4,m)})(t)= t^{m}M(\overline{r})(t)+ S(t)= t^{m}(t^6-t^5-2t)+ (-2t^5-t+1),$$ 
and Proposition~\ref{prop_PF}(2) says that $\lim_{m \to \infty} \lambda(\beta_{(4,m)}) \approx 1.45109$. 
\end{ex}

We are now ready to show  Proposition~\ref{prop_decreasing-example}. 
\medskip

{\it Proof of Proposition~\ref{prop_decreasing-example}.} 
Recall the tree $\mathcal{Q}_{(m_1, \cdots, m_{k+1})}$ and the tree map $q_{(m_1, \cdots, m_{k+1})}$ 
used in the proof of Proposition~\ref{prop_pseudo-Anosov}. 
For $i$ even, $\mathcal{Q}_{(m_1, \cdots, m_{k+1})}$ is also the combined tree obtained from the triple 
$(\mathcal{Q}_{(m_1, \cdots, m_{i-1})}, \mathcal{G}_{m_i,-}, \mathcal{Q}_{(m_{i+1}, \cdots, m_{k+1})})$. 
Then, $q_{(m_1, \cdots, m_{k+1})}$ is also the combined tree map given by 
$$\widehat{q}_{(m_1, \cdots, m_{i-1})}  g_{m_i} \widehat{q}_{(m_{i+1}, \cdots, m_{k+1})}: 
\mathcal{Q}_{(m_1, \cdots, m_{k+1})} \rightarrow \mathcal{Q}_{(m_1, \cdots, m_{k+1})},$$
where  $\widehat{q}_{(m_1, \cdots, m_{i-1})}$ and $\widehat{q}_{(m_{i+1}, \cdots, m_{k+1})}$ 
are suitable extensions of $q_{(m_1, \cdots, m_{i-1})} $ and  $q_{(m_{i+1}, \cdots, m_{k+1})}$ respectively. 
By Proposition~\ref{prop_decreasing}, the claim holds.

For $i$ odd, $\mathcal{Q}_{(m_1, \cdots, m_{k+1})}$ is the combined tree obtained from 
($\mathcal{Q}_{(m_1, \cdots, m_{i-1})}, \mathcal{G}_{m_i,+}, \mathcal{Q}'_{(m_{i+1}, \cdots, m_{k+1})})$, 
where  $\mathcal{Q}'_{(m_{i+1}, \cdots, m_{k+1})}$ is the tree obtained from $\mathcal{Q}_{(m_{i+1}, \cdots, m_{k+1})}$ by the horizontal reflection. 
Then, the proof  is similar to that for the even case.  
$\Box$
\medskip

We turn to the proof of  Theorems~\ref{thm_recursive-example} and \ref{thm_lim-example}.
\medskip

{\it Proof of Theorem~\ref{thm_recursive-example}.} 
By Proposition~\ref{prop_BH}  and by the proof of Proposition~\ref{prop_pseudo-Anosov}, 
the dilatation of $\beta_{(m_1, \cdots, m_{k+1})}$ is the largest root of $M(q_{(m_1, \cdots, m_{k+1})})(t)$. 
Fixing  $m_1  ,\cdots, m_k \ge 1$, 
let ${\mathcal R}_{(m_1, \cdots, m_k)}$ and 
$\overline{r}_{(m_1, \cdots, m_k)}$ be the dominant tree and the dominant tree map for $\{q_{(m_1, \cdots, m_{k+1})}\}_{m_{k+1} \ge 1}$, 
and we set 
$$R_{(m_1, \cdots, m_k)}(t)= M(\overline{r}_{(m_1, \cdots, m_k)})(t).$$
By Proposition~\ref{prop_recipe}, 
\begin{eqnarray}
\label{equation_key1}
M(q_{(m_1, \cdots, m_{k+1})})(t)&=& t^{m_{k+1}} R_{(m_1, \cdots, m_k)}(t)+ S_{(m_1, \cdots, m_k)}(t)\ \mbox{and}
\\
\label{equation_key2}
{M(q_{(m_1, \cdots, m_{k+1})})}_*(t)&=& t^{m_k+1}U_{(m_1, \cdots, m_k)}(t) + {R_{(m_1, \cdots, m_k)}}_*(t),
\end{eqnarray}
where $S_{(m_1, \cdots, m_k)}(t)$ is the recessive polynomial for $\{M(q_{(m_1, \cdots, m_{k+1})})(t)\}_{m_{k+1} \ge 1}$ 
and $U_{(m_1, \cdots, m_k)}(t)$ is the dominant polynomial for $\{{M(q_{(m_1, \cdots, m_{k+1})})}_*(t)\}_{m_{k+1} \ge 1}$.

\begin{claim}
\label{claim_key1}
We have 
\medskip
\\
(1) $S_{(m_1, \cdots, m_k)}(t) =  (-1)^{k+1} {R_{(m_1, \cdots, m_k)}}_*(t)$ and 
\medskip
\\
(2) $U_{(m_1,  \cdots, m_k)}(t)= (-1)^{k+1} R_{(m_1, \cdots, m_k)}(t)$. 
\end{claim}

{\it Proof.}  
It is enough to show Claim~\ref{claim_key1}(1). 
For  if (1) holds,  by (\ref{equation_key1}) we have 
$$M(q_{(m_1, \cdots, m_{k+1})})(t)= (-1)^{k+1} M(q_{(m_1, \cdots, m_{k+1})})_*(t).$$ 
This together with (\ref{equation_key1}), (\ref{equation_key2}) implies Claim~\ref{claim_key1}(2). 

We prove  Claim~\ref{claim_key1}(1) by an induction on $k$. 
For $k=1$, this holds \cite[Theorem~3.20(1)]{HK}. 
We assume Claim~\ref{claim_key1}(1) up to $k-1$. 
Then, we have $S_{(m_1, \cdots, m_{k-1})}(t) =  (-1)^{k} {R_{(m_1, \cdots, m_{k-1})}}_*(t)$, 
\begin{eqnarray}
\label{equation_Q}
M(q_{(m_1, \cdots, m_k)})(t)= (-1)^k {M(q_{(m_1, \cdots, m_k)})}_*(t) \ \mbox{and}\ 
U_{(m_1, \cdots, m_{k-1})}(t) =  (-1)^{k} R_{(m_1, \cdots, m_{k-1})}(t). 
\end{eqnarray}
For $k \ge 2$, 
the transition matrix $M(q_{(m_1,\cdots, m_{k+1})})$ has the block form: 
\bigskip

$M(q_{(m_1,\cdots, m_{k+1})})=$
\begin{scriptsize}
$\bordermatrix{
              &1& & n_1+1                                                &                                                  & n_2+1&\cdots &n_{k-1}+1 & & n_k+1& & & n_{k+1} \cr
              &  & &                                                             &                                                 &             &           &                   & &             & & &               \cr
              &  & &                                                             &M(q_{(m_1,\cdots, m_{k})})  &             &           &                   & &             & & &               \cr
              &  & &                                                            &                                                   &             &           &                   & & 1            & & &               \cr
              &  & &                                                            &                                                    &             &           &                              & & 1            & & &               \cr 
n_k+1  &1   & & 2                                                      &                                                     &2           &\cdots &  2              & &  1           &1 & &               \cr
              &   & &                                                           &                                                     &              &            &                 & &                & &\ddots &               \cr
              &   & &                                                           &                                                      &              &            &                 & &                & &            & 1          \cr
n_{k+1}        &1   & & 2                                                      &  &2           &\cdots &  2              & &  2           & & &               \cr              
},$
\end{scriptsize}
\bigskip
\\
where $n_j= m_1 + \cdots + m_j + j$ for $1 \le j \le k+1$. 
Note that the last edge of the tree ${\mathcal R}_{(m_1, \cdots, m_j)}$ is numbered $n_j$. 
In this case the polynomial $S_{(m_1, \cdots, m_k)}(t)$ is the determinant of a matrix: 
\bigskip

\begin{scriptsize}
$$\left(
\begin{array}{ccccccccc|cccc}
                   &   &        &                                                                        & & & &         & &                  & & & \\
                    &  &        &tI- M(\overline{r}_{(m_1, \cdots, m_{k-1})}) & & & &         &  &                   & & & \\
                    &  &        &                                                                        & & & &         & & -1           & & & \\ \hline
                     &  &       &                                                                        & & & &           & & t&\ddots & & \\
                    &   &       &                                                                         & & & &           & & &\ddots         &-1 & -1\\
                - 1&   &-2  &                                                                         & -2& &\cdots & &-2   &           & &t &-1 \\
                  -1&  &-2  &                                                                         &-2 & &\cdots &  & -2  &           &   &    &-2 \\ 
\end{array}
\right). 
$$
\end{scriptsize}
\bigskip
\\
Subtract the second last row from the last row of this matrix, and let $A= (a_{i,j})$ be the resulting matrix. 
Applying the determinant expansion with respect to the last row of $A$, we have
\begin{eqnarray*}
 S_{(m_1, \cdots, m_k)}(t) = |A|
 = \sum_{j=1}^{n_k+1}(-1)^{n_k+1+j} a_{n_k+1,j}|A_{n_{k+1},j}|
 \end{eqnarray*}
 where $A_{i,j}$ is the matrix obtained by $A$ with row $i$ and column $j$ removed.  
 Since $a_{n_k+1,j}=0$ for $j \ne n_k, n_k+1$ and $a_{n_k+1,n_k}=-t$, $a_{n_k+1,n_k+1}=-1$, 
 we have 
 \begin{eqnarray*}
 |A| =  t |A_{n_k+1,n_k}|- |A_{n_k+1,n_k+1}|. 
 \end{eqnarray*}
 We note that  $|A_{n_k+1,n_k+1}|= M(q_{(m_1, \cdots, m_k)})(t)$. 
 For the computation of $|A_{n_k+1,n_k}|$, 
 subtract the second last row of $A_{n_{k+1},n_{k+1}-1}$ from the last row, and 
for the resulting matrix, apply the determinant expansion of the last column successively. 
Then, we obtain 
$$|A_{n_k+1,n_k}|= -t^{m_k-1} R_{(m_1, \cdots, m_{k-1})}(t)+ S_{(m_1, \cdots, m_{k-1})}(t).$$ 
Thus, 
$$S_{(m_1, \cdots, m_k)}(t)= -M(q_{(m_1, \cdots, m_k)})(t)- t^{m_k} R_{(m_1, \cdots, m_{k-1})}(t)+ t S_{(m_1, \cdots, m_{k-1})}(t).$$
In the same manner we have 
$${R_{(m_1, \cdots, m_k)}}_*(t)= {M(q_{(m_1, \cdots, m_k)})}_*(t)+ t^{m_k} U_{(m_1, \cdots, m_{k-1})}(t)-t{R_{(m_1, \cdots, m_{k-1})}}_*(t). $$
By using (\ref{equation_Q}), these two equalities imply Claim~\ref{claim_key1}(1).  
This completes the proof. 
\medskip

We now turn to proving  Theorem~\ref{thm_recursive-example}. 
We will prove an inductive formula for $R_{(m_1, \cdots, m_k)}(t)$. 
It is not hard to show that $R_{(m_1)}(t)= t^{m_1+1}(t-1)-2t$.

For $k \ge 2$, one can  verify 
\begin{equation}
\label{equation_recursive}
R_{(m_1, \cdots, m_k)}(t)= t M(q_{(m_1, \cdots, m_k)})(t)- t^{m_k} R_{(m_1, \cdots, m_{k-1})}(t)+t S_{(m_1, \cdots, m_{k-1})}(t). 
\end{equation}
Substitute the two equalities 
\begin{eqnarray*}
M(q_{(m_1, \cdots, m_k)})(t)&=& t^{m_k} R_{(m_1, \cdots, m_{k-1})}(t)+ (-1)^k{R_{(m_1, \cdots, m_{k-1})}}_*(t)\ \mbox{and}
\\
S_{(m_1, \cdots, m_{k-1})}(t)&=& (-1)^k{R_{(m_1, \cdots, m_{k-1})}}_*(t)
\end{eqnarray*}
into (\ref{equation_recursive}), then we find the inductive formula 
$$R_{(m_1, \cdots,m_k)}(t) = t^{m_k}(t-1) R_{(m_1, \cdots,m_{k-1})}(t)+ (-1)^k 2t {R_{(m_1, \cdots,m_{k-1})}}_*(t).$$
This completes the proof of Theorem~\ref{thm_recursive-example}. 
$\Box$
\medskip

{\it Proof of Theorem~\ref{thm_lim-example}.} 
To begin with, we show 

\begin{claim}
\label{claim_easy}
We have
\medskip
\\
(1) $\displaystyle\lim_{m_i \to \infty} \lim_{m_{i+1} \to \infty} \cdots \lim_{m_{k+1} \to \infty} \lambda(\beta_{(m_1,  \cdots, m_{k+1})})=1$ for $i = 1$ and 
\medskip
\\
(2) $\displaystyle\lim_{m_i \to \infty} \lim_{m_{i+1} \to \infty} \cdots \lim_{m_{k+1} \to \infty} \lambda(\beta_{(m_1,  \cdots, m_{k+1})})=
\lambda(R_{(m_1, \cdots, m_{i-1})}(t)) >1$ for $i \ge 2$. 
\end{claim}

{\it Proof.} 
(1) 
By Theorem~\ref{thm_recursive-example} and by Lemma~\ref{lem_asymptotic-root1} 
$$\displaystyle\lim_{m_{k+1} \to \infty} \lambda(\beta_{(m_1, \cdots, m_{k+1})}) = 
\lambda (R_{(m_1, \cdots, m_k)}(t)).$$
Recall that $R_{(m_1, \cdots, m_i)}(t)= M(\overline{r}_{(m_1, \cdots, m_i)})(t)$. 
It is not hard to see that the matrix  $M(\overline{r}_{(m_1, \cdots, m_i)})$ for $i \ge 1$ is PF, and hence 
the largest root of $R_{(m_1, \cdots, m_i)}(t)$ is greater than $1$. 
Then, by using the inductive formula of $R_{(m_1, \cdots, m_i)}(t)$ in Theorem~\ref{thm_recursive-example} 
together with Lemma~\ref{lem_asymptotic-root1}, we have  
$$\displaystyle\lim_{m_i \to \infty} \lambda(R_{(m_1, \cdots, m_i)}(t)) = \lambda(R_{(m_1, \cdots, m_{i-1})}(t))$$
for $i \ge 2$. 
Since $R_{(m_1)}(t) = t^{m_1+1}(t-1)-2t$, we have $\displaystyle\lim_{m_1 \to \infty} \lambda(R_{(m_1)}(t)) = 1$ by  Lemma~\ref{lem_asymptotic-root2}. 
This completes the proof of (1). 

(2) The proof is identical to  that of (1). 
This completes the proof of Claim~\ref{claim_easy}. 
\medskip

Claim~\ref{claim_easy}(1) says that for any $\lambda>1$ there exists an integer $m_i(\lambda)$ for each $i$ with $1 \le i \le k+1$ 
such that $\lambda(\beta_{m_1(\lambda), \cdots, m_{k+1}(\lambda)})< \lambda$. 
Set $m= \max\{m_i(\lambda)\ |\ i = 1, \cdots k+1\}$. 
By  Proposition~\ref{prop_decreasing-example} $\lambda(\beta_{(m_1, \cdots, m_{k+1})})< \lambda$ whenever $m_i > m$.  
This completes the proof of Theorem~\ref{thm_lim-example}(1). 

The proof of Theorem~\ref{thm_lim-example}(2) is identical to that of (1), but using Claim~\ref{claim_easy}(2) instead of Claim~\ref{claim_easy}(1). 
$\Box$
\medskip

We show the existence of two kinds of families of pA mapping classes with arbitrarily small dilatation and with arbitrarily large volume. 
\medskip

{\it Proof of Proposition~\ref{prop_trivial-construction}.} 
There exists a family of pseudo-Anosov mapping classes $\psi_n$ of $\mathcal{M}(D_n)$ such that  
$$\lim_{n \to \infty} \lambda(\psi_n) = 1.$$ 
It suffices to show that for any pA mapping class  $\phi \in \mathcal{M}(\Sigma_{g,p})$, 
there exists a family of pA mapping classes $\hat{\phi}_n \in \mathcal{M}(\Sigma_{g,p(n)})$ 
such that  the dilatation of $\hat{\phi}_n$ is same as $\phi$ and the volume of $\hat{\phi}_n$ goes to $\infty$ as $n$ goes to $\infty$.

 Let $\Phi \in \phi$ be a pA homeomorphism. 
 Since the set of periodic orbits of $\Phi$ is dense on $\Sigma_{g,p}$, one can find a periodic orbit of $\Phi$, say $Q= \{q_1, \cdots, q_{s}\}$. 
 Now puncture each point of  $Q$, then the pA mapping class $\phi' \in \mathcal{M}(\Sigma_{g,p+s})$ induced 
 by $\phi$ satisfies $\lambda(\phi')= \lambda(\phi)$. 
 On the other hand, $\mathrm{vol}(\phi') > \mathrm{vol}(\phi)$ 
 since ${\Bbb T}(\phi)$ is a complete hyperbolic manifold obtained topologically by filling a cusp of ${\Bbb T}(\phi')$ with a solid torus \cite[Section~6]{Thu2}. 
  The volume of any cusp is bounded below uniformly. 
  Thus, if we puncture periodic orbits of $\Phi \in \phi$ successively, we obtain a family of pA mapping class with the desired property.   
  $\Box$
\medskip

\noindent
 Finally, we show Theorem~\ref{thm_dil-vol}. 
 \medskip
 
{\it Proof of Theorem~\ref{thm_dil-vol}}. 
By Theorem~\ref{thm_lim-example}(1), for each integer $k \ge 1$, the dilatation of $\beta_{(m_1,\cdots, m_{k+1})}$ goes to $1$ as 
$m_1, \cdots, m_{k+1}$ all go to $\infty$. 
Thus, it suffices to show that the volume of $\beta_{(m_1,\cdots, m_{k+1})}$ goes to $\infty$ as $k$ goes to $\infty$. 

One verifies that  $\overline{\beta_{(m_1, \cdots, m_{k+1})}}$ is a $2$ bridge link as in Figure~\ref{fig_general_2bridge}. 
In particular it is an alternating link with  twist number  $k+1$.  
Theorem~1  in \cite{Lac} tells us that for each $m_1, \cdots, m_{k+1} \ge 1$, 
$$\mathrm{vol}(\beta_{(m_1,\cdots, m_{k+1})}) > \frac{1}{2} (k-1) v_3, $$ 
where $v_3$ is the  volume of a regular ideal tetrahedron. 
This completes the proof.  
$\Box$

\noindent
\hspace{10cm}
Eiko Kin and Mitsuhiko Takasawa
\\
\hspace{10cm}
Department of Mathematical and 
\\
\hspace{10cm}
Computing Sciences
\\
\hspace{10cm}
Tokyo Institute of Technology 
\\
\hspace{10cm}
Tokyo,  Japan
\\
\hspace{10cm}
 kin@is.titech.ac.jp 
\\
\hspace{10cm}
takasawa@is.titech.ac.jp

\end{document}